\documentclass[12pt]{article}
\usepackage{amsmath}
\usepackage{amsthm}
\usepackage{amssymb}
\def\pf{{\bf Proof. }}
\def\beq{\begin{equation}}
\def\eeq{\end{equation}}

\def\Si{\Sigma}

\def\g{\gamma}

\def\e{\epsilon}
\def\o{\omega}
\def\O{\Omega}

\def\al{\alpha}
\def\bt{\beta}
\def\ro{\rho}

\def\t{\tilde}

\def\nd{\noindent}
\def\d{\delta}
\def\i{\iota}
\def\p{\partial}
\def\hf{\hfill{$\Box$}}
\def\<{\leq}
\def\>{\geq}

\def\sub{\subset}
\newtheorem{thm}{Theorem}[section]
\newtheorem{lem}{Lemma}[section]
\newtheorem{prop}{Proposition}[section]
\newtheorem{cor}{Corollary}[section]

\newtheorem{rem}{Remark}[section]

\begin{document}
\title{\bf The Degree of Symmetry of Certain Compact Smooth Manifolds II}
\author{\sc Bin Xu}
\date{}
\maketitle
\begin{abstract}
\noindent In this paper, we give the sharp estimates for the degree of symmetry
and the semi-simple degree of symmetry of certain four dimensional fiber
bundles by virtue of the rigidity theorem of harmonic maps due to Schoen and
Yau. As a corollary of this estimate, we compute the degree of symmetry and the
semi-simple degree of symmetry of ${\Bbb C}P^2\times V$, where $V$ is closed
smooth manifold admitting a real analytic Riemannian metric of non-positive
curvature. In addition, by the Albanese map, we obtain the sharp estimate of
the degree of symmetry of a compact smooth manifold with some restrictions on
its one dimensional cohomology.
\end{abstract}

\noindent
{\bf Mathematics Subject Classification:} Primary 57S15; Secondary 53C44.\\

\noindent
 {\bf Key  Words and Phrases:} Degree of symmetry, fiber bundle,
non-positive curvature, harmonic map, first cohomology.

\section{Introduction}
\ \ \ \
Let $M^n$ be
a closed, connected and smooth
$n$-manifold and $N(M^n)$ the {\it degree of
symmetry} of $M^n$, that is, the maximum of the dimensions of the
isometry groups of all possible Riemannian metrics on $M^n$.
Of course, $N(M)$ is the maximum of the dimensions of the compact Lie groups
which can act effectively and smoothly on $M$.
The {\it semi-simple degree of symmetry} $N_s(M)$ is defined similarly, where
we consider only actions of semi-simple  compact Lie groups on $M$.
The following is well known:
\begin{equation}
\label{equ:est}
N(M^n)\leq n(n+1)/2.
\end{equation}
In addition, if the equality holds, then $M^n$ is diffeomorphic to the
standard sphere $S^n$ or the real projective space ${\Bbb R}P^n$. In
\cite{KM} H. T. Ku, L. N. Mann, J. L. Sicks and J. C. Su obtained similar
results on a product manifold $M^n=M^{n_1}_1\times M^{n_2}_2\ (n\geq 19)$
where $M_i$ is a compact connected smooth manifold of dimension
$n_i$: they showed that
\begin{equation}
\label{equ:pro}
N(M)\leq n_1(n_1+1)/2+n_2(n_2+1)/2,
\end{equation}
and that if the equality holds,
then $M^n$ is a product of two spheres, two real projective spaces or
a sphere and a real projective space.
A preliminary lemma for the proof of Ku-Mann-Sicks-Su's results
claims that if $M^n$ $(n\geq 19)$ is a compact connected smooth $n$-manifold
which is not diffeomorphic to the complex projective space ${\Bbb C}P^m$
$(n=2m)$, then
\begin{equation}
\label{equ:be}
N(M^n)\leq k(k+1)/2+(n-k)(n-k+1)/2
\end{equation}
holds for each $k\in {\bf N}$ such that the $k$-th Betti number $b_k$ of $M$ is nonzero.

Let $V$ be a closed, connected and smooth  manifold
which admits a real analytic Riemannian metric of non-positive curvature.
It was noted in Remark 1.2 of \cite{Xu} that by the results in \cite{CR3} and
\cite{LY} the following holds:\\

\nd{\bf Fact N} $N(V)$ equals the rank of Center$\,\pi_1(V)$ and any
connected compact Lie groups which can act effectively and smoothly on $V$ must be a torus group. \\

\nd Let $E$ be a closed smooth fiber
bundle over $V$ with connected fiber $F$. In Theorem 1.1 of \cite{Xu} the
author generalized partially Ku-Mann-Sicks-Su's result \eqref{equ:pro} by
showing the corresponding sharp estimates of $N(E)$ and $N_s(E)$ by assuming
that the fiber $F$ satisfies various topological properties.
In particular, part of the statements of Theorem 1.1 and Corollary 1.1 in
\cite{Xu} says: \\

\nd {\bf Fact F} Suppose that $E$ is oriented and that $F$ is an orientable $4m$-manifold
($m\geq 5$) of nonzero signature. Then the following statements hold{\rm :}
\begin{equation}
\label{equ:sig}
N(E)\leq N(V)+4m(m+1),\ \ N_s(E)\leq 4m(m+1)\ .
\end{equation}
In particular, if $V$ is orientable, then
\begin{equation}
\label{equ:cpm}
N({\Bbb C}P^{2m}\times V)=N(V)+4m(m+1),\ \ N_s({\Bbb C}P^{2m}\times V)=4m(m+1)\ .
\end{equation}
\\

\nd The case of $1\leq m\leq 4$ could not be covered in \cite{Xu}
because the author used Ku-Mann-Sick-Su's result \eqref{equ:be}, in which
the dimension of the manifold is assumed to be $\geq 19$. In this paper we
shall show that Fact F also holds for $m=1$. That is,

\begin{thm}
\label{thm:sig} Let $V$ be a closed, connected and smooth manifold
which can admit a real analytic Riemannian metric of non-positive
curvature and $E$ be a closed smooth fiber bundle over $V$ such
that the fiber $F$ of $E$ is connected. Suppose that $E$ is
oriented and that the fiber $F$ has dimension $4$ and has nonzero
signature. Then the followings hold:
\begin{equation}
\label{equ:sig4}
N(E)\leq N(V)+8,\ \ N_s(E)\leq 8\ .
\end{equation}
In particular, if $V$ is oriented, then
\begin{equation}
\label{equ:cp2}
N({\Bbb C}P^2\times V)=N({\Bbb C}P^2)+N(V)=N(V)+8,\ \ N_s({\Bbb C}P^2\times V)=
N_s({\Bbb C}P^2)=8\ .
\end{equation}
\end{thm}

In fact,  the assumption in \eqref{equ:cp2} that $V$
is oriented can be removed by the following

\begin{thm}
\label{thm:uno}
Let $V$ be a closed, connected and smooth manifold which can admit a real analytic
Riemannian metric of non-positive curvature and $E$ be a closed smooth fiber
bundle over $V$ such that the fiber $F$ of $E$ is connected.
Suppose that the fiber $F$ has dimension $4$ and is not cobordant mod 2 to
either 0 or ${\Bbb R}P^4$. Then both \eqref{equ:sig4} and \eqref{equ:cp2} hold.
\end{thm}

\begin{rem}
{\rm
The assumption  in Theorem \ref{thm:sig} that $F^4$ is orientable and has nonzero signature
in Theorem \ref{thm:sig} is independent of the assumption in Theorem \ref{thm:uno} that $F^4$ is not cobordant mod 2 to
either 0 or ${\Bbb R}P^4$ in Theorem \ref{thm:uno}. For examples, the
oriented $4$-manifold ${\Bbb C}P^2\sharp {\Bbb C}P^2$ has signature 2 and is cobordant mod 2 to 0;
${\Bbb R}P^2\times {\Bbb R}P^2$ is a non-orientable
$4$-manifold which is not cobordant to either $0$ or ${\Bbb R}P^4$. However,
since the group $\O_4$ of the oriented cobordism class of 4-manifolds is
isomorphic to ${\Bbb Z}$ and the isomorphism is given by the signature, an
oriented compact 4-manifold having zero signature is cobordant to zero in $\O_4$.

By Remark 1.4 in \cite{Xu} the connectedness of $F$ is necessary for the
validity of Theorems \ref{thm:sig} and \ref{thm:uno}.
}
\end{rem}

D. Burghelea and R. Schultz \cite{BS} showed that $N_s(M)=0$ if there exist
$\al_1,\cdots,\al_n$ in $H^1(M;{\Bbb R})$ with $\al_1\cup\cdots\al_n\not=0$.
In Theorem 1.2 of \cite{Xu} Burghelea-Schultz's result was generalized to the
following\\

\nd{\bf Fact C} Let $M$ be an $n$-dimensional closed connected smooth
manifold. If there exist $\al_1,\cdots,\al_k$ in $H^1(M;{\Bbb R})$ with
$\al_1\cup\cdots\cup\al_k\not=0$,
then the followings hold:
\[N(M)\leq (n-k+1)(n-k)/2+k\ ,\]
\[ N_s(M)
\left\{
\begin{array}{rl}
\leq (n-k+1)(n-k)/2 & {\rm if}\ n-k>1\\
=0 & {\rm otherwise}\ .
\end{array}
\right.
\]
\\

Further assuming $b_1(M)>k$, we obtain the following

\begin{thm}
\label{thm:cup} Let $M$ be an $n$-dimensional closed connected
smooth manifold and $k\geq 3$ be an integer. If the first Betti
number $b_1(M)$ of $M$ is greater than $k$ and there exist
$\al_1,\cdots,\al_k$ in $H^1(M;{\Bbb R})$ with
$\al_1\cup\cdots\cup\al_k\not=0$ in $H^k(M;{\Bbb R})$, then
\begin{equation}
\label{equ:btk}
N(M)\leq
\left\{
\begin{array}{rl}
(n-k+1)(n-k)/2+k-2 & {\rm if}\ n\geq k+3\\
n & {\rm if}\ n=k+2 \ {\rm or}\ k+1\\
n-2 & {\rm if}\ n=k\ .
\end{array}
\right.
\end{equation}
\end{thm}

\begin{rem}
{\rm It is implied by the assumption of Theorem \ref{thm:cup} that
$M$ has dimension $\geq k$. We make the assumption of $k\geq 3$ in
Theorem \ref{thm:cup} because of some known facts in \cite{Xu}.
Precisely speaking, the statements (ii) and (iii) of Theorem 1.2
in \cite{Xu} give best possible estimates for Case $k=0$ or $1$
and Case $k=2$, respectively. Moreover, to get those estimates
there, the author only need to assume that $b_1(M)>k$ because the
non-vanishing property of the cup product
$\al_1\cup\cdots\cup\al_k$ is superfluous in the sense that it
does not give any more restrictions to the degree of symmetry of
$M$.}
\end{rem}

\begin{rem}
{\rm By the definition of degree of symmetry, it is easy to see that for
a product manifold $M_1\times M_2$, where $M_i$ is a compact connected smooth
manifold, the following holds:
\begin{equation}
\label{equ:pd}
N(M_1\times M_2)\geq N(M_1)+N(M_2)\ .
\end{equation}
Let $\Si_g$ be the oriented closed surface of genus $g$ and
$M^n=S^{n-k}\times T^{k-2}\times \Si_g$ $(n\geq k+3,\ g\geq 2$). Then $M^n$
satisfies the assumption of Theorem \ref{thm:cup}. Since by Fact N
$N(T^{k-2}\times \Si_g)=k-2$, by \eqref{equ:pd} and Theorem \ref{thm:cup} we
obtain the equality
\[N(S^{n-k}\times T^{k-2}\times \Si_g)=(n-k+1)(n-k)/2+k-2\ ,\]
combining which with the equalities
\[N(T^{n})=n,\ \ N(T^{n-2}\times \Si_g)=n-2,\]
we can see that the estimate
\eqref{equ:btk} is best possible.
}
\end{rem}

This paper is organized as follows. In Section \ref{sec:pre}, we prepare for
the following sections. In particular, we cite some results in \cite{Xu} and
prove a key lemma (cf Lemma \ref{lem:dim4}) for Theorems
\ref{thm:sig} and \ref{thm:uno}. In Section \ref{sec:sig}, we prove Theorem
\ref{thm:sig} (\ref{thm:uno}) with the help of this key lemma and
the oriented (unoriented) cobordism theory. In Section \ref{sec:cup},
 we prove Theorem \ref{thm:cup} by virtue of the unique
continuation property of harmonic maps.


\section{Preliminaries}
\label{sec:pre}
\ \ \ For a closed Riemannian manifold $M$ let $I^0(M)$ be the identity
component of the isometry group of $M$.
The following proposition will provide the frame for the proof of Theorems
\ref{thm:sig}, \ref{thm:uno}.

\begin{prop}
\label{prop:mod} {\rm (cf \cite{SY} Theorem 4 )} Suppose that $M,\
N$ are closed real analytic Riemannian manifolds and the sectional
curvature of $N$ is non-positive. Suppose that $h:M\to N $ is a
surjective harmonic map and its induced map $h_*:\pi_1(M)\to
\pi_1(N)$ is surjective. Then the space of surjective harmonic
maps homotopic to $h$ is represented by $\{\bt\circ h| \bt\in
I^0(N)\}$, where $I^0(N)$ is a torus group of dimension equaling
both the rank of {\rm Center}$\,\pi_1(N)$ and the degree of
symmetry of $N$.
\end{prop}

We cite a topological result from \cite{Xu}.

\begin{prop}
\label{prop:fib}
{\rm (cf \cite{Xu} Proposition 3.1)}
Let $p_0:E\to B$ be a fiber bundle over a compact connected
smooth manifold $B$ such that the fiber of $E$ is also connected.
Then any continuous map homotopic to $p_0:E\to B$ is surjective.
\end{prop}

We cite a lemma in \cite{Xu}, which is also necessary for Theorems
\ref{thm:sig}, \ref{thm:uno}.

\begin{lem}
\label{lem:eff} {\rm (cf \cite{Xu} Lemma 2.1)} Let $M$ be a
connected Riemannian manifold and $f$ a smooth map from $M$ to a
smooth manifold $N$. Suppose that $y\in N$ is a regular value of
$f$ and $F$ is a connected component of the submanifold
$f^{-1}(y)$ of $M$. If an isometry $\al$ of $M$ satisfies that
$h\circ\al=h$ and that $\al(x)=x$  for any $x\in F$, then $\al$ is
the identity map of $M$.
\end{lem}

\begin{lem}
\label{lem:dim4}
{\rm (Key lemma of Theorems \ref{thm:sig} and \ref{thm:uno})}
Let $Y$ be a closed connected smooth {\rm 4}-manifold not diffeomorphic to
either $S^4$ or ${\Bbb R}P^4$. Then $N(Y)\leq 8$.
The equality $N(Y)=8$ holds if and only if $Y$ is diffeomorphic to
${\Bbb C}P^2$. Moreover, $N_s({\Bbb C}P^2)=8$.
\end{lem}
\nd\pf By \eqref{equ:est} $N(Y)\leq 9$. Then $N(Y)\leq 8$ follows
from Theorem A$'$ in Ishihara \cite{Is} which claims that there
exists no 4-dimensional Riemannian manifold having a
$9$-dimensional isometry group.  If $Y$ is a Riemannian manifold
whose isometry group has dimension 8, then by Theorem 5 in
Ishihara \cite{Is} $Y$ is a K{\"a}hlerian space with positive
constant holomorphic sectional curvatures. Since the holomorphic
sectional curvature of a K{\" a}hler manifold determines
completely its Riemannian curvature tensor (cf \cite{Zhe} Lemma
7.19.), $Y$ has positive sectional curvature and then by the
theorem of Synge (cf \cite{CE} Theorem 5.9.)  $Y$ is simply
connected. By the theorem of Cartan-Ambrose-Hicks (cf \cite{CE}
Theorem 1.36), $Y$ is isometric to the K{\" a}hler manifold ${\Bbb
C}P^2$ with the Fubini-Study metric. Since the compact Lie group
SU(3) acting isometrically on ${\Bbb C}P^2$ is semi-simple,
$N_s({\Bbb C}P^2)=N({\Bbb C}P^2)=8$. \hf
\\

We do some preparations for the proof of Theorem \ref{thm:cup} in the
following.

For a compact oriented Riemannian manifold $M$ with nonzero first Betti number $b_1(M)$,
let ${\cal H}$ be the real vector space of all harmonic 1-forms on $M$ and $\nu$ the natural
projecion from the universal covering ${\t M}$ of $M$. For $x_0\in {\t M}$, set $p_0=\nu(x_0)$.
We define a smooth map
${\t a}:{\t M}\to {\cal H}^*$ from ${\t M}$ to the dual space  ${\cal H}^*$
of  ${\cal H}$ by a line integral
\[{\t a}(x)(\o)=\int_{x_0}^{x}{\nu^*\o}.\]
For $\sigma\in \pi_1(M)$
\[{\t a}(\sigma x)={\t a}(x)+\psi(\sigma)\]
holds, where $\psi(\sigma)(\o)=\int_{x_0}^{\sigma x_0}{\nu^* \o}$,
so that $\psi$ is a homomorphism from $\pi_1(M)$ into ${\cal H}^*$
as an additive group. It is a fact that $\Delta =\psi(\pi_1(M))$
is a lattice in the vector space ${\cal H}^*$, and clearly this
vector space has a natural Euclidean metric from the global inner
product of forms on $M$. With the quotient metric, we call the
torus $A(M)={\cal H}^*/\Delta$ the $Albanese\ torus$ of the
Riemannian manifold $M$. By the above relation between ${\t a}$
and $\psi$, we obtain a map $a:M\to A(M)$ satisfying ${\t a}(x)\in
a\circ \nu(x)$ for any $x\in {\t M}$. We call the map $a$ the
$Albanese\ map$. From the very construction of $a$, we see that
the map it induces on fundamental groups
\[a_*:\pi_1(M)\to\pi_1(A(M))\]
is surjective and that $a^*$ maps the space of harmonic 1-forms on $A(M)$ isomorphically onto
${\cal H}$. By Corollary 1 in \cite{NS}, the Albanese map is harmonic.
Set
\[r_a:=\max\{{\rm rank}\,da(p)|p\in M\}\ .\]

\begin{lem}
\label{lem:albrk}
{\rm (cf \cite{Xu} Lemma 4.3)}
Let $M$ be an n-dimensional oriented compact Riemannian manifold
with nonzero first Betti number $b_1$. Let $a:M\to A(M)$ be the Albanese map.
Suppose there exist $\al_1,{\cdots},\al_k$ in $H^1(M;{\Bbb R})$
with $\al_1\cup{\cdots}\cup\al_k\not=0$ in $H^k(M;{\Bbb R})$.
Then $r_a\geq k$ holds.\\
\end{lem}

\begin{lem}
\label{lem:rel}
{\rm (cf \cite{Xu} Lemmata 4.1 and 4.2)}
Let $M$ be a non-orientable compact manifold and
$\pi:M'\to M$ be its orientable double covering. Then the following statements hold
{\rm:}\\
{\rm (i)}\ \ \ $N(M)\leq N(M')$\\
{\rm (ii)}\ \ $b_1(M)\leq b_1(M')$\\
{\rm (iii)}\ If $M$ has the property that there exist $k$ one dimensional
real cohomology classes $\al_1,\cdots,\al_k$ of $M$ such that
$\al_1\cup\cdots\cup\al_k$ is nonzero in $H^k(M;{\Bbb R})$, then
so does $M'$.
\end{lem}


\section{Proof of Theorems \ref{thm:sig} and \ref{thm:uno}}
\label{sec:sig} \nd{\sc Proof of Theorem \ref{thm:sig}}\ \ For the
proof of \eqref{equ:sig4}, by Corollary \ref{cor:rdg} in Appendix
we have only to show that for any real analytic Riemannian metric
on $E$ and any compact semi-simple subgroup $G$ of $I^0(E)$, the
following inequalities hold:
\begin{equation}
\label{equ:G}
\dim I^0(E)\leq N(V)+8,\ \dim G\leq 8 \ .
\end{equation}
Since the fiber $F$ is connected, the fiber bundle projection $p:E\to V$
induces a surjective map $p_*:\pi_1(E)\to \pi_1(V)$. Using a well known result
by Eells-Sampson \cite{ES}, we see that there exist harmonic maps homotopic
to $p:E\to V$. By Proposition \ref{prop:fib}, each of them is surjective
and then satisfies the assumptions of Proposition \ref{prop:mod}. Taking a
harmonic map $h:E\to V$ homotopic to $p:E\to V$, by Proposition \ref{prop:mod},
for any $\al\in I^0(E)$ we can find a unique $\ro(\al)\in I^0(V)$ with
$h\circ\al=\ro(\al)\circ h$. We see that $\ro: I^0(E)\to I^0(V)$ is a Lie group
homomorphism. Since $G$ is semi-simple and $I^0(N)$ is a torus group, the restriction of $\ro$ to $G$
must be trivial. That is, $G$ is contained in Ker$\,\ro$. Therefore, the proof of \eqref{equ:G}
is completed if we can show that Ker $\ro$, which acts isometrically on $E$,
has dimension $\leq 8$.

Choosing a smooth homotopy $P: E\times [0,\ 1]\to V$ between $p$
and $h$, we can also choose a regular value $y$ of $P$ by Sard's
theorem since $P$ is surjective. Then we have the following\\

\nd {\bf Claim 1} $P^{-1}(y)$ is a oriented submanifold in $E\times [0,\,1]$
with boundary $p^{-1}(y)+h^{-1}(y)$. That is, there exists an oriented
cobordism in $E$ between $F$ and $h^{-1}(y)$. \\

\nd{\it Proof of Claim 1}\ \ Since $y$ is the regular value of $P$ and
 $P^{-1}(y)$ is non-empty, it is easy to see
that $P^{-1}(y)$ is a submanifold of $E\times [0,\,1]$ with boundary
\[\p P^{-1}(y)\cap E\times 0+\p P^{-1}(y)\cap E\times 1=p^{-1}(y)+h^{-1}(y)
\cong F+h^{-1}(y)
 .\]
Note that the normal bundle of $P^{-1}(y)$ in $E$ is trivial.
Since $E$ is orientable, so is $P^{-1}(y)$. We proved the claim.
\\

By Hirzebruch's signature theorem (cf \cite{Hz} Theorem 8.2.2) and Claim 1,
up to the sign difference, the signature of ${\t F}$ equals that of $F$
so that there exists a connected component $F^*$ of ${\t F}$ having nonzero
signature. Hence $F^*$ is not diffeomorphic to either $S^4$ or
${\Bbb R}P^4$. By Lemma \ref{lem:eff} Ker $\ro$ acts effectively on
$F^*$. Moreover, Lemma \ref{lem:dim4} tells us that Ker\ $\ro$ has dimension $\leq 8$.
We complete the proof of \eqref{equ:G} and \eqref{equ:sig4}. Since ${\Bbb C}P^2$ has signature 1,
\eqref{equ:cp2} follows from \eqref{equ:sig4} and \eqref{equ:pd}.
\hf\\

\nd{\sc Proof of Theorem \ref{thm:uno}} \ \ Repeating the part of
the proof of Theorem \ref{thm:sig} before Claim 1, we can see that
$P^{-1}(y)$ is a submanifold of $E\times [0,\,1]$ with boundary
$p^{-1}(y)+h^{-1}(y)$. That is, $h^{-1}(y)$ is cobordant mod 2 to
$F$. Since $F$ is not cobordant mod 2 to either 0 or ${\Bbb
R}P^4$, there exists a connected component $F^*$ of $h^{-1}(y)$
such that $F^*$ is not diffeomorphic to either $S^4$ or ${\Bbb
R}P^4$. Since Ker $\ro$ acts effectively on $F^*$, by Lemma
\ref{lem:dim4} Ker $\ro$ has dimension $\leq 8$. Therefore the
inequalities in \eqref{equ:sig4} hold. To show \eqref{equ:cp2}, we
only need to show ${\Bbb C}P^2$ is not cobordant mod 2 to zero or
${\Bbb R}P^4$, which follows from that $w_2^2[{\Bbb R}P^4]=0$ and
$w_2^2[{\Bbb C}P^2]\not=0$.
\hf\\


\section{Proof of Theorem \ref{thm:cup}}
\label{sec:cup} \nd{\sc Proof of Theorem \ref{thm:cup}}\ \ By
Lemma \ref{lem:rel}, we may assume $M$ is an oriented Riemannian
manifold. Let  $a:M\to A(M)$ be the Albanese map and $b_1$ the
first Betti numer of $M$. By Corollary \ref{cor:rdg} we have only
to consider the analytic Riemannian metric on $M$. For any $\g\in
I^0(M)$, $a\circ\g$ is also a harmonic map from $M$ to the
Albanese torus $A(M)$ and homotopic to $a$. By Lemma 3 in
\cite{NS} there is a unique translation $\ro(\g)$ of the torus
$A(M)$ such that
\[a\circ\g=\ro(\g)\circ a\ .\]
Then we have a Lie group homomorphism $\ro: I^0(M)\to T^{b_1}$, where the torus
group $T^{b_1}$ is the translation group of the Albanese torus $A(M)$.

Remember $r_a=\max\{{\rm rank}\,da(p)|p\in M\}$. We claim that
\begin{equation}
\label{equ:ker}
\dim\ {\rm Ker}\ \ro\leq \frac{1}{2}(n-r_a+1)(n-r_a),\ \ \dim\ {\rm Im}\ \ro\leq r_a.
\end{equation}
\nd{\it Proof of \eqref{equ:ker}}\quad Since the proof in Lemma
2.3 \cite{Xu} has some ambiguity, we give a clear and rigorous
proof here. We first prove $\dim\,{\rm Im}\,\ro\leq r_a$. Let $r$
be the dimension of Im$\,\ro$. As a connected subgroup of
$T^{b_1}$, Im$\,\ro$ is an $r$-dimensional torus group acting
freely on $a(M)$ by the definition of $\ro$. Choose a point $x$ in
$M$. Then $y:=a(x)$ is in $a(M)$ and the orbit $\bigl({\rm
Im}\,\ro\bigr)(y)$ of $y$ with respect to the ${\rm Im}\,\ro$
action is an $r$-dimensional subtorus contained in $a(M)$. More
precisely, by the definition of $\ro$, we have
\[\bigl({\rm Im}\,\ro\bigr)(y)=a\Bigl(I^0(M) (x)\Bigr),\]
where $I^0(M)(x)$ is the orbit of $x$ with respect to the $I^0(M)$ action on $M$.
Let $a_1$ be the restriction of the Albanese map $a$ to the submanifold $I^0(M) (x)$ of $M$. Then
$r_a\geq r_{a_1}=\dim\,{\rm Im}\,\ro=r$.

Then we show that $K:={\rm Ker}\,\ro$ has dimension $\leq
(n-r_a+1)(n-r_a)/2$ by investigating the $K$ action on $M$. Since
$M$ is connected, the orbit space $M/K$ with the induced topology
is also connected. By Theorem 4.27 \cite{Ka}, the union set of all
principal orbits of the $K$ action is an open dense subset of $M$.
Taking an open set $U\subset M$, at any point of which the map
$da$ has rank equal to $r_a$, we can choose a point $p$ in $U$
such that the orbit $K(p)$ is principal. By the definition of
$\ro$, $K(p)$ is contained in the inverse image of $a(p)$ under
$a$. Since ${\rm rank}\,da(p)=r_a$, the submanifold $K(p)$ has
dimensional not exceeding $n-r_a$. Since $K$ acts effectively on
the principal orbit $K(p)$, the dimensional of $K$ can not exceed
$(n-r_a+1)(n-r_a)/2$ by \eqref{equ:est}.
By now we have completed the proof of \eqref{equ:ker}.\\

We see from Lemma \ref{lem:albrk} that $r_a\geq k$.

\nd{\it Case 1}\ \ If $r_a\geq k+1$,
then from \eqref{equ:ker}
\[\dim\ I(M)=\dim\ {\rm Ker}\ \ro+\dim\ {\rm Im}\ \ro
\leq \frac{1}{2}(n-k-1)(n-k)+k+1.\]

\nd{\it Case 2}\ \ Suppose $r_a=k$ in what follows. We claim that
dim Im $\ro$ will be less than $k-1$, which together with
\eqref{equ:ker} imply the following estimate
\[\dim\ I^0(M)\leq\frac{1}{2}(n-k+1)(n-k)+k-2. \]
Otherwise, suppose $\dim\ {\rm Im}\ \ro\geq k-1$. Remember that
the Lie group ${\rm Im}\ \ro$ acting on $A(M)$ in fact acts on the
image $a(M)$ of $a$. Hence we can assume that there exists a
subgroup $T^{k-1}$ of the translation group $T^{b_1}$ which acts
freely and isometrically on $a(M)$. Since both $M$ and $A(M)$ are
real analytic, a theorem of Morrey \cite{Mo} shows that the
harmonic mapping $a$ is in fact real analytic. By well-known
theorems in real analytic geometry \cite{Lo} we know that both $M$
and $A(M)$ can be triangulated so that $a(M)$ is a compact
connected simplicial subcomplex of dimension $k$ in $A(M)$. We
write the orbit space of the free and isometric $T^{k-1}$ actions
on $A(M)$ and $a(M)$ by $A(M)/T^{k-1}$ and $a(M)/T^{k-1}$
respectively, in which the former is in fact also a flat torus of
dimension $b_1-k+1$. Since the natural projection map $\pi:A(M)\to
A(M)/T^{k-1}$ is totally geodesic, we see that by a result in
\cite{EL} the composition map $\pi\circ a: M\to A(M)/T^{k-1}$ is a
harmonic map, whose image is $a(M)/T^{k-1}$, the orbit space of
the free $T^{k-1}$ action on the simplicial subcomplex $a(M)$ of
dimension $k$ in $A(M)$. Hence $a(M)/T^{k-1}$, the image of
$\pi\circ a$ in $A(M)/T^{k-1}$, has dimension $1$ so that the
differential of harmonic map $\pi\circ a$ has rank $\leq 1$ at any
point of $M$. By the unique continuation property of the harmonic
maps (cf \cite{Sam} Theorem 3 ), we see that $\pi\circ a$ maps $M$
onto a closed geodesic of $A(M)/T^{k-1}$, which means that $a(M)$
is a principal $T^{k-1}$-bundle over $S^1$. Since $S^1$ is
connected, there exists a section on this bundle so that $a(M)$ is
a trivial $T^{k-1}$-bundle, i.e. a $k$-dimensional torus. This
contradicts the surjectivity of the homomorphism $a_*:\pi_1(M)\to
\pi_1(A(M))\cong {\Bbb Z}^{b_1}$ ($b_1>k$). Hence, we proved the
claim. In particular, if the dimension $n$ of the manifold $M$
equals $k$, then $r_a$ must be $k$ and the estimate
$\dim\,I^0(M)\leq k-2$ follows from the claim.

Combining Cases 1 and 2, we can see that the dimension of $I^0(M)$
is dominated by the maximum of those two integers
\[(n-k-1)(n-k)/2+k+1, \quad\quad (n-k)(n-k+1)/2+k-2,\]
provided $\dim\,M=n\geq k+1$. By some simple computations, we
complete the proof.
\hf\\

\appendix
\section{Real Analytic Group Action}
In this appendix, we will prove the following theorem.
\begin{thm}
\label{thm:rdg}
Let $G$ be a compact Lie group acting smoothly and effectively on a compact
smooth manifold $M$. Then there exists a real analytic manifold $M'$ on which
$G$ acts real analytically such that there exists a $G$-isomorphism
between $(M,G)$ and $(M', G)$. That is, there is a diffeomorphism $f:M\to M'$
which satisfies for any $x\in M$ and $g\in G$
\[f(gx)=g\bigl(f(x)\bigr)\ .\]
\end{thm}

Although the theorem should not be new, the author provides a proof here since he
does not know any reference of the theorem.
We put the proof of Theorem \ref{thm:rdg} afterward. Firstly, from it
we have the following
\begin{cor}
\label{cor:rdg}
The degree of symmetry of $M$ equals the maximum of the dimensions
of the isometry groups of all the real analytic Riemannian metrics on $M$.
\end{cor}
\nd \pf
Let $G$ be a compact Lie group acting smoothly and effectively on a compact
smooth manifold $M$. By Theorem \ref{thm:rdg} there exists another
real analytical $G$ action on the unique real analytic structure of $M$
compatible to the existed smooth structure o $M$. Moreover, the new $G$ action
is equivariant to the old one on the smooth structure of $M$.
Thus the new one is also effective. Taking a real analytic Riemannian metric
on $M$, by the invariant integration on $G$ we can construct a new real
analytic one on which $G$ acts isometrically. \hf\\

Let $f:M\to N$ be a smooth map between smooth manifolds $M$ and $N$.
It is {\it transverse} to a submanifold $A\sub N$ if and only if
whenever $f(x)=y\in A$, then the tangent space to $N$ at $y$ is spanned by
the tangent space to $A$ at $y$ and the image of the tangent space to $M$ at
$x$. That is,
\[ T_yA+df(T_xM)=T_yN .\]
\begin{lem}
{\rm (cf Theorem 1.3.3  \cite{Hir})}
\label{lem:tsv}
Let $f:M\to N$ be a smooth map and $A\sub N$ a submanifold of codimension $\ell$.
If $f$ is transverse to $A$, then $f^{-1}(A)$ is a submanifold of $M$ of
codimension $\ell$.
\end{lem}

\begin{lem}
{\rm (cf Theorem 4.12  \cite{Ka})}
\label{lem:ee}
Let $G$ be a compact Lie group and $M$ a compact manifold on which $G$ acts
smoothly. Then there exists a representation space $(V,\mu)$ of $G$ and a
smooth $G$-embedding $\iota:M\to V$. That is, for any $x\in M$ and $g\in G$,
\[\i(gx)=\mu(g)\bigl(\i(x)\bigr)\ .\]
Moreover, if the $G$-action on $M$ is effective, then the
representation $(V,\mu)$ of $G$ is faithful.
\end{lem}

Let $G$ be a compact group acting smoothly on two manifolds $M$ and $N$.
A smooth map $f:M\to N$ is a {\it G-map} if and only if for any
$x\in M$ and $g\in G$ the following holds:
\[g\bigl(f(x)\bigr)=f(gx)\ .\]
\begin{lem}
{\rm (An equivariant version of Theorem 2.5.2 \cite{Hir})}
\label{lem:eap}
Let $G$ be a compact Lie group acting isometrically on Euclidean spaces
${\Bbb R}^q$ and ${\Bbb R}^s$. Let $M\sub {\Bbb R}^q$ be a $G$-invariant compact
submanifold of codimension $>0$ and $E$ a $G$-invariant tubular neighborhood
of $M$ in ${\Bbb R}^q$. Let $f:E\to W$ be a smooth $G$-map into a $G$-invariant
open set $W$ of ${\Bbb R}^s$. Let $v:{\Bbb R}^q\to {\Bbb R}$ be a smooth
$G$-invariant function with support in $E$, equal to $1$ on a $G$-invariant
compact
neighborhood $K$ of $M$. Set $h(x)=v(x)f(x)=v(x)(f_1(x),\cdots,f_s(x))$.
Let $\d: {\Bbb R}^q\to {\Bbb R},\,\d(x)=\exp\,(-|x|^2)$.
Let $C=1/\int_{{\Bbb R}^q} \delta$.
Let $\e>0$. Then for $k>0$ sufficiently large,
\[\psi(x)=(\psi_1(x), \cdots,\psi_s(x)):=(h_1(x),\cdots,h_s(x))*(Ck^q\d(kx))\]
is an analytic $G$-map and satisfies $||\psi-f||_{C^1,\,K}<\e$.
\end{lem}
\nd\pf The proof of the analytic property of $\psi$ and the estimate
$||\psi-f||_{C^1,\,K}<\e$ for $k>0$ large enough is straightforward.
We have only to show that $h*\d$ is $G$-invariant. Since $v(gx)=v(x)$,
$h=vf:E\to {\Bbb R}^s$ is $G$-map. For any $g\in G$, since $g$ acts on
${\Bbb R}^q$ isometrically and $\d$ is a radial function on ${\Bbb R}^q$,
\begin{eqnarray*}
(h*\d)(gx)&=&\int_{{\Bbb R}^q} h(y)\d(gx-y)\,dy
=\int_{{\Bbb R}^q} h(gz)\d(g(x-z))\,dz\\
&=&\int_{{\Bbb R}^q} h(gz)\d(x-z)\,dz=
g\biggl(\int_{{\Bbb R}^q} h(z)\d(x-z)\,dz\biggr)\\
&=&g\Bigl(h*\d (x)\Bigr)\ .
\end{eqnarray*}
\hf
\\

\nd{\sc Proof of Theorem \ref{thm:rdg}}\quad By Lemma \ref{lem:ee}, there
exists a faithful representation space $V$ of $G$ and a
$G$-embedding $\i:M\to V$. By the invariant integration on $G$,
we can induce a $G$-invariant inner product $(\, ,\,)$ on $V$, equipped
with which $V$ becomes a Euclidean space ${\Bbb R}^q$ and $G$ becomes a
subgroup of O\,($q$). We shall not distinguish $M$ and $\iota(M)$ in what follows.
Let $k$ be the codimension of $M$ in ${\Bbb R}^q$.
Since $M$ is a $G$-invariant submanifold, by Theorem 4.8 in \cite{Ka}
there exists a $G$-invariant normal tubular neighborhood $E$ of $M$ in ${\Bbb R}^q$
which can be identified with a $G$-invariant
neighborhood of the zero section of the normal
bundle of $M$ in ${\Bbb R}^q$. Let $p:E\to M$ be the restriction of the bundle projection to $E$,
which is a $G$-map.

Let $G_{q,k}$ be the Grassmann manifold of $k$-dimensional linear subspaces
of ${\Bbb R}^q$ and $E_{q,k}\to G_{q,k}$ be the Grassmann bundle, the fiber
of $E_{q,k}$ over the $k$-plane $P\sub {\Bbb R}^q$ is the set of pairs
$(P,x)$ where $x\in P$. Then the $G$-action on ${\Bbb R}^q$ induces the natural
real analytic actions on $G_{q,k}$ and $E_{q,k}$ respectively
such that the bundle projection $E_{q,k}\to G_{q,k}$ is a $G$-map.
Let $h:M\to G_{q,k}$ be the map sending $x\in M$ to the $k$-plane normal
to $M$ at $x$ and $f:E\to E_{q,k}$ be the natural
map covering $h$; thus
\[f(y)=(h\circ p(y),y)\in E_{q,k}\sub G_{q,k}\times {\Bbb R}^q \ .\]
Since $G$ acts isometrically on $M$, as a linear map
on ${\Bbb R}^q$, $dg=g$ maps the $k$-plane normal to $M$ at
$x$ to the one normal to $M$ at $gx$. Therefore, both $h$ and $f$ are $G$-maps.
Moreover,
$f$ is transverse to the zero section $G_{q,k}\sub E_{q,k}$ and
\[f^{-1}(G_{q,k})=M\ .\]

Now we embed $E_{q,k}$ analytically in ${\Bbb R}^s$ with $s=q^2+q$. For this
it suffices to embed $G_{q,k}$ in ${\Bbb R}^{q^2}$. This is done by mapping
a $k$-plane $P\in G_{q,k}$ to the linear map ${\Bbb R}^q\to{\Bbb R}^q$ given
by the orthogonal projection on $P$. There exists a natural isometric
$G$-action on ${\Bbb R}^s$ such that the embedding $E_{q,k}\sub {\Bbb R}^s$ is a $G$-map. Then we can find a $G$-invariant normal tubular
neighborhood $W$ of $E_{q,k}$. Let $\Pi:W\to E_{q,k}$ be the real
analytic $G$-invariant projection.

Let $f':E\to W$ be the extension of the map $f:E\to E_{q,k}$ to
$W$. It follows from Lemma \ref{lem:eap} that $f'$ can be approximated
near $M$ by an analytic $G$-map $\psi:E\to W$. Then $\phi=\Pi\circ \psi$ is
an analytic $G$-invariant approximation of $f:E\to E_{q,k}$.
Put $M'=\phi^{-1}(G_{q,k})$. If $\phi$ is
sufficiently $C^1$ close to $f$, then $\phi$ is also transverse to
$G_{q,k}$, so by Lemma \ref{lem:tsv} $M'$ is a real analytic submanifold of codimension $k$ in $E\subset {\Bbb R}^q$.
Moreover, the restriction of $G$-map $p:E\to M$ to $M'$ is a $G$-isomorphism from $M$ to $M'$. \hf\\

\nd
{\bf Acknowledgements}

\nd This study was supported in part by the Japanese Government
Scholarship, the JSPS Postdoctoral Fellowship for Foreign
Researchers, the FRG Postdoctoral Fellowship and the Program of
Visiting Scholars at Nankai Institute of Mathematics. The author
would like to express his gratitude to Professor Mikiya Masuda for
his valuable comments in Remark 1.1. The author also thank the
referee, many valuable and concrete comments from whom made this
paper more readable.

\noindent {\small \sc Nankai Institute of Mathematics, Nankai
University, Tianjin 300071 China\\
Former address: Department of Mathematics, Johns Hopkins University, Baltimore MD 21218}\\
{\small \it E-mail address}: \tt{bxu@math.jhu.edu}

\end{document}